\documentclass{procDDs}
\newcommand{\ben}{\begin{eqnarray}}
\newcommand{\een}{\end{eqnarray}}
\newcommand{\la}{\label}

\title{Novel representation of the general Heun's functions. Back to the beginning}

\author{\textbf{Plamen~P~Fiziev}}
{Sofia University Foundation for Theoretical and Computational Physics and Astrophysics, Boulevard
5 James Bourchier, Sofia 1164, Bulgaria\\
and\\
BLTF, JINR, 141980 Dubna, Moscow region, Russia}                 
{fizev@phys.uni-sofia.bg\,\,\,and\,\,\,
fizev@theor.jinr.ru}                                   
\begin {document}

\maketitle

\index{Fiziev, P.P.}                              
\begin{abstract}
We study a novel type of solutions of the general Heun's equation,
based on its symmetric form.
We derive the symmetry group of this equation which is a proper extension of the Mobius group.
The new series solution treat simultaneously and on an equal footing all singular points.
\end{abstract}

\section{Introduction}

The popular form of the general Heun's equation
\begin{multline}
H''+\left({\frac {\gamma_{{}_G}} {z}}+{\frac{\delta_{{}_G}}{z-1}}+{\frac{\epsilon_{{}_G}}{z-a_{{}_G}}}\right)H'+\\
{\frac{\alpha_{{}_G}\beta_{{}_G} z-\lambda}{z(z-1)(z-a_{{}_G})}}H = 0,\,\,\,\,\,
\gamma_{{}_G}\!+\delta_{{}_G}\!+\epsilon_{{}_G}\!=\!\alpha_{{}_G}\!+\!\beta_{{}_G}\!+\!1;\nonumber
\label{dHeunG}
\end{multline}
treats asymmetrically its regular singular points $0, 1, a_{{}_G}, \infty$
with indices \cite{Heun,Ron,SL}
\begin{multline}
\hskip 1.truecm \{ 0, 1-\gamma_{{}_G}\},\quad \{ 0, 1-\delta_{{}_G} \},\\
\{ 0, \gamma_{{}_G}+\delta_{{}_G}-\alpha_{{}_G}-\beta_{{}_G} \},\quad
\{ \alpha_{{}_G}, \beta_{{}_G} \}.
\end{multline}
This causes difficulties in numerical calculations outside
the circle $|z|<1$, $z \in \mathbb{C}$, especially for $|z|\gg1$.

The symmetric form of the general Fuchsian equation \cite{Klein,Forsyth} when applied to
the special case of four regular singular points $z_{j=1...4}\in \mathbb{C}$
with the indices $\{\alpha_j,\beta_j\}:{j=0,...,4}$, $\sum_{j=1}^4 \left(\alpha_j+\beta_j\right)=2$ is
\begin{multline}
\mathcal{W}''+\left(\sum_{j=1}^4{\frac{1-\alpha_j-\beta_j}{z-z_j}}\right)\mathcal{W}'+\\
{\frac 1 {P(z)}\left(\lambda+\sum_{j=1}^4{\frac{q_j}{z-z_j}}\right)}\mathcal{W} = 0,
\la{dWN}
\end{multline}
\begin{multline}
P(z)=\prod_{j=1}^4(z-z_j)=\sum_{n=0}^4(-1)^n\sigma_{{}_{4-n}} z^n,\\
q_j=\alpha_j\beta_j P^\prime(z_j)\,\,\, \text{for}\,\,\,{j=1...4}.
\end{multline}

Making use of the transformation
\ben
\mathcal{W}(z) = \mathcal{F}(z)\prod_{j=1}^4 \left(z-z_j\right)^{\nu_j}, \,\, \sum_{j=1}^4 \nu_j=0,
\la{nu_transf}
\een
and imposing {\it symmetric} constraints on the indices
\ben
\alpha_j+\beta_j={\tfrac 1 2}, \,\, \alpha_j\beta_j P^\prime(z_j)=q_j,\,\, \text{for}\quad j=1..4,
\la{index}
\een
we obtain the following simple symmetric form of the general Heun's equation
\begin{multline}
\mathcal{F}^{\prime\prime}+{\frac 1 2}\left(\sum_{j=1}^4{\frac 1{z-z_j}}\right)\mathcal{F}^\prime+\\
{\frac 1 {P(z)}\left(\lambda+\sum_{j=1}^4{\frac{q_j}{z-z_j}}\right)}\mathcal{F} = 0,
\la{dF}
\end{multline}
or in a self-adjoint form:
\ben
\left(P(z)\right)^{1/2}\left(\left(P(z)\right)^{1/2}\mathcal{F}^\prime\right)^\prime +
\left(\lambda+Q(z)\right)\mathcal{F} = 0,
\label{dFsa}
\een
where
\begin{multline}
Q(z)=\sum_{j=1}^4{\frac{q_j}{z-z_j}}=\\
-{\tfrac 1 {16 P(z)}}
\sum_{j=1}^4\! \left(\sin(2\chi_j)\right)^2\partial_z P(z_j)\partial_{z_j}P(z)
\end{multline}

Here we use new uniformization parameters $\chi_{j=1..4} \in \mathbb{C}$:
\begin{multline}
\hskip 1.truecm \alpha_j={\tfrac 1 2} \left(\cos\chi_j\right)^2,\quad \beta_j={\tfrac 1 2} \left(\sin\chi_j\right)^2,\nonumber\\
\quad q_j={\tfrac 1 {16}}\left(\sin(2\chi_j)\right)^2 \prod_{k\neq i}^4\left(z_j-z_k\right).\hskip 1.truecm
\la{alphabeta}
\end{multline}

\section{The orthogonality and the invariance group of solutions}
\hskip ,37truecm {\bf Proposition 1:}
For any two solutions $\mathcal{F}_{\lambda_{1,2}}(z)$  of the Eq. \eqref{dFsa} with  $\lambda_1\neq\lambda_2$ we have
\begin{multline}
\int_{\mathcal{L}_{ij}}\mathcal{F}_{\lambda_1}(z)\mathcal{F}_{\lambda_2}(z)d\mu(z)=0,\\
d\mu(z)=\left(P(z)\right)^{-1/2}dz.\hskip .35truecm
\end{multline}

Here $\mathcal{L}_{ij}\in\mathbb{C}$ is any contour which starts at the singular point $z_j$ and ends at the singular point $z_j$
without going through the other singular points $z_{k\neq i,j}$. Besides, the singular boundary conditions
\begin{multline}
\left(P(z)\right)^{1/2}
\left(\mathcal{F}_{\lambda_1}(z)\mathcal{F}_{\lambda_2}^\prime(z)-
\mathcal{F}_{\lambda_2}(z)\mathcal{F}_{\lambda_1}^\prime(z)\right)\,\upharpoonleft_{z_j,z_j}\quad\\=0
\label{boundary}
\end{multline}
are supposed to be satisfied. $ \blacktriangleleft$

As seen, in our symmetric formulation of the general Heun's equation \eqref{dFsa}
the parameter $\lambda$ plays natural role of eigenvalue in the
corresponding two-singular-point boundary problems for a self-adjoined operator
\cite{Ron,SL,KLS,LS}. This is also consistent with the profound results of \cite{Smirnov}.

{\bf Proposition 2:}
Equation \eqref{dF} is invariant under the {\it extension} $\widehat{\mathfrak{G}}_{Mobius}$
of the standard Mobius group $\mathfrak{G}_{Mobius}$. This extension
acts on the functions of $10$ variables $\mathcal{F}\left(z;z_1,...,z_4;q_1,...,q_4;\lambda\right)$
and is produced by the following basic transformations:

\begin{description}
  \item[(i)] Complex translations with  $\zeta \in \mathbb{C}$:
\begin{multline}
z\to  z+\zeta; \quad z_j\to z_j+\zeta:\, j=1,...,4;\\
q_j\to q_j:\, j=1,...,4; \quad
\lambda \to \lambda.\hskip  1.1truecm
\label{translationN4}
\end{multline}
  \item[(ii)] Complex dilatations with  $t\in \mathbb{C}$, $t\neq 0$:
\begin{multline}
z\to t\, z; \quad z_j\to t\, z_j:\, j=1,...,4; \\
q_j\to t^3 q_j:\, j=1,...,4;\quad \lambda \to t^2 \lambda.\hskip  .55truecm
\label{rescalingN4}
\end{multline}
  \item[(iii)] Inversion
\begin{multline}
\hskip .truecm z\to 1/ z; \quad z_j\to 1/ z_j:\, j=1,...,4; \\
q_j \to -q_j/\left( z_j^2 \sigma_4 \right) :\, j=1,...,4;\hskip  1.truecm \\
\lambda \to \Big(\lambda -\sum_{j=1}^4 q_j/z_j\Big)/\sigma_4.\hskip  .6truecm \,  \blacktriangleleft
\label{inversionN4}
\end{multline}
\end{description}

Using proper compositions of these basic transformations (see \cite{Smirnov,AF} for
${\mathfrak{G}}_{Mobius}$ ) we are able to construct
a representation of the whole extended Mobius group $\widehat{\mathfrak{G}}_{Mobius}$ that acts on the solutions
$\mathcal{F}\left(z;z_1,...,z_4;q_1,...,q_4;\lambda\right)$ of Eq. \eqref{dF}
without bringing us outside of the variety of these solutions.
Hence, the extended group $\widehat{\mathfrak{G}}_{Mobius}$
is the group of invariance of the variety of solutions to Eq. \eqref{dF}.

\section{The Taylor series expansion of the solutions}
Our next step is to adopt the following basic assumption which is of crucial importance for further work:
\ben
z_{j=0,...,4}\neq 0.
\la{zero}
\een
Then the function $\mathcal{F}(z)$
is an analytical one in the vicinity of the point $z=0$ and has an absolutely convergent
Taylor series expansion
\begin{multline}
\mathcal{F}(z)\equiv\mathcal{F}(z;z_1,...z_4;q_1,...q_4;\lambda) =\\
\sum_{n=0}^\infty f_n(z_1,...z_4;q_1,...q_4;\lambda) z^n
\label{taylorF}
\end{multline}
with the coefficients $f_n(z_1,...z_4;q_1,...q_4;\lambda)$ defined by the nine-term recurrence relation
\ben
f_n+\sum_{k=1}^8 r_{n-k}f_{n-k}=0,
\label{frecurrence}
\een
which can be obtained from Eqs. \eqref{dF} and \eqref{taylorF}.

We shall take advantage of the freedom to put the singular points $z_{j=0,...,4}$ in the proper places in the complex plane $\mathbb{C}$
for simplifying, as much as possible, the coefficients in recurrence \eqref{frecurrence} and thus, the very solutions $\mathcal{F}(z)$.
This can be done in a symmetric way by imposing additional conditions on the elementary symmetric functions $\sigma_{j=1,2,3,4}$.
Using the proper Mobius transformation one is able to impose three independent constraints on $z_{j=0,...,4}$
without changing the problem.
Obvious simple choice is to reduce the quartic equation $P(z)=0$ to the following biquadratic one: $z^4- 2\cos(2\phi)z^2+1=0$
with the roots
\ben
z_1=e^{i\phi},\,\,z_2=-e^{-i\phi},\,\,z_3=-e^{i\phi},\,\,z_4=e^{-i\phi},
\la{z_j}
\een
by imposing three symmetric constraints
\ben
\sigma_1=\sigma_3=0,\,\, \sigma_4=1,
\la{B_sigma_constraits}
\een
and replacing $\sigma_2= -2\cos(2\phi)$ with one more complex uniformization parameter $\phi\in \mathbb{C}$.

The meaning of the new variable $\phi$ is revealed by the formula for the invariant $a(z_1,z_2,z_3,z_4)$
of the Mobius transformation -- the so called {\it cross-ratio}.
In our problem it acquires the form
\begin{multline}
a={\frac {(z_1-z_3)(z_2-z_4)}{(z_2-z_3)(z_1-z_4)}}={\frac 1{\left(\sin\phi\right)^2}}
\quad \Rightarrow\\
\sigma_2=-2\left(1-{\frac 2 a}\right).
\label{a}
\end{multline}
As a result, one obtains for the coefficients in recurrence \eqref{frecurrence}
in the following form:
\begin{align}
r_{n-1}&=0,\label{r134:a}\\
r_{n-2}&={\tfrac 1 {n(n-1)}}\Big(\lambda- {\tfrac i 4}\,\sin(2\phi)\,\rho_2\Big)+\nonumber\\
&4\left(1-{\tfrac  5 n}+{\tfrac  3 2}{\tfrac  1 {n-1}}\right)\cos(2\phi), \label{r134:b}\\
r_{n-3}&=-{\tfrac 1 {n(n-1)}}{\tfrac i 4}\,\sin(2\phi)\,\rho_3, \label{r134:c}\\
r_{n-4}&={\tfrac 1 {n(n-1)}}\Big(\!-2\lambda \cos(2\phi)+ {\tfrac i 4}\,\sin(2\phi)\,\rho_4\Big)+\nonumber\\
&2\left(1-{\tfrac  {16} n}+{\tfrac  9 {n-1}}\right)\left((\cos(2\phi)^2+1\right),\label{r134:d}\\
r_{n-5}&={\tfrac 1 {n(n-1)}}{\tfrac i 4}\,\sin(2\phi)\,\rho_5, \label{r134:e}\\
r_{n-6}&={\tfrac \lambda {n(n-1)}}-\nonumber\\
&4\left(1-{\tfrac  {33} n}+{\tfrac  {45} 2}{\tfrac  1 {n-1}}\right)\cos(2\phi),\label{r134:f}\\
r_{n-7}&=0,\label{r134:g}\\
r_{n-8}&=1-{\tfrac  {56} n}+{\tfrac  {42} {n-1}}.\label{r134:h}
\end{align}
where we introduce the following four functions of five variables $\{\phi,\chi_1,\chi_2\chi_3,\chi_4\}$:
\begin{align}
\rho_2=\,&\left(\left(\sin(2\chi_1)\right)^2+\left(\sin(2\chi_3)\right)^2\right)-\nonumber\\
&\left(\left(\sin(2\chi_2)\right)^2+\left(\sin(2\chi_4)\right)^2\right),\nonumber\\
\rho_3=\,&e^{-i\phi}\left(\left(\sin(2\chi_1)\right)^2-\left(\sin(2\chi_3)\right)^2\right)+\nonumber\\
&e^{i\phi}\,\left(\left(\sin(2\chi_2)\right)^2-\left(\sin(2\chi_4)\right)^2\right),\nonumber\\
\rho_4=\,&e^{2i\phi}\,\left(\left(\sin(2\chi_1)\right)^2+\left(\sin(2\chi_3)\right)^2\right)-\nonumber\\
&e^{-2i\phi}\left(\left(\sin(2\chi_2)\right)^2+\left(\sin(2\chi_4)\right)^2\right),\nonumber\\
\rho_5=\,&e^{i\phi}\,\,\left(\left(\sin(2\chi_1)\right)^2-\left(\sin(2\chi_3)\right)^2\right)+\nonumber\\
&e^{-i\phi}\left(\left(\sin(2\chi_2)\right)^2-\left(\sin(2\chi_4)\right)^2\right).\nonumber
\end{align}

\section{The circular case}
There exist two quite different cases of positions of the singular points $z_{j=0,...,4}$.

The first (general) case is the one in which the singular points $z_{j=0,...,4}$ do not lie
on any circle in the complex plane. We call it  {\it the non circular case}.
In the non circular case, the theory of solutions of Eq. \eqref{dF} is more complicated.
It is not enough developed for the popular form of the general Heun's equation
with complex $a_{{}_G} \notin \mathbb{R}$, see the Introduction.

The second (special) case is the one when all regular singular
points $z_{j=0,...,4}$ of Eq. \eqref{dF} lie
on some circle $\mathfrak{C}\in \mathbb{C}$.
We will call this case {\it the circular case}.
A particular and natural such case is the one when $z_{j=0,...,4}\in \mathbb{R}$, i.e.,
all singular points are real.

In the circular case, one is able to move all singular points $z_{j=0,...,4}$
on any other circle $\mathfrak{\tilde C}\in \mathbb{C}$ using Mobius transformation
which preserves the circular property
and transforms any circle in the complex plane into another circle,
see \cite{Smirnov,AF}. Especially, in \cite{Smirnov}
the circular case for the general Heun's equation with $a_{{}_G} \in \mathbb{R}$ was substantially
elaborated using the standard choice \eqref{dHeunG}
and without any relation with the choice \eqref{z_j} in Eq. \eqref{dF}.

The value of the invariant cross-ratio \eqref{a} $a\in\mathbb{R}$ is real
for any four complex points $z_{j=1,2,3,4}$ on a circle in $\mathbb{C}$.
Then, from relation \eqref{a} follows that in the circular case the angle $\phi\in\mathbb{R}$ is real and
the singular points \eqref{z_j} lie on the unit circle with the center $z=0$.

Thus, according to the basic results of the standard analytic theory of ordinary differential equations
\cite{Golubev,CL}, we obtain our key result:

{\bf Proposition 3:} In the circular case, the series \eqref{taylorF} with coefficients \eqref{r134:a}-\eqref{r134:h}
and $\phi\in\mathbb{R}$ are absolutely convergent inside the unit circle, i.e., for any $z\in \mathbb{C}$ with $|z|<1$.\,$\blacktriangleleft$

{\bf Corollary:} In the circular case, the four regular singular points $z_{j=1,2,3,4}$
of the general Heun's equation \eqref{dF} can be treated simultaneously
and equally from inside the unit circle using the Taylor series \eqref{taylorF}.

Next important step is to restrict Proposition 4 (iii) to the circular case.

{\bf Proposition 4:}
In the circular case, equation  \eqref{dF} preserves its form if one makes the following substitutions
\ben
z&\to& 1/z,\quad z_j \to 1/z_j, \\
F(z)&\to& F(1/z),\quad P(z)\to P(1/z),\nonumber\\
\lambda &\to& \lambda -{\tfrac i 4}\sin(2\phi)\rho_2,\quad q_j \to -q_j/z_j^2.\nonumber
\la{inverssion_circular}
\een
This way we obtain from solutions \eqref{taylorF} new solutions of Eq. \eqref{dF} in the form of the
Laurent series expansions which are absolutely convergent for any $z\in \mathbb{C}$ with $|z|>1$.\,$\blacktriangleleft$

{\bf Corollary:} In the circular case, the four regular singular points $z_{j=1,2,3,4}$
of the general Heun's equation \eqref{dF} are mapped under inversion on the same points,
removed to the initial positions $z_{j=4,3,2,1}$, respectively. Hence,
one can treat all singular points simultaneously and equally from outside the unit circle
using the corresponding  Laurent series, described in Proposition 6.

As a final result, in the circular case we reach a totally symmetric treatment of the singular points
$z_{j=1,2,3,4}$  in the whole complex plane $\mathbb{\tilde C}$.

\section{Concluding remarks}
In the present talk, we introduced and studied a novel representation of the general Heun's functions.
It is based on the symmetric form of the Heun's differential equation
yielded by a further development of the Felix Klein
symmetric form of the Fuchsian equations with an
arbitrary number $N\geq 4$ of regular singular points.

The basic relations for the general Heun's equation with $N=4$ are derived and discussed in detail.

We derived the symmetry group of this equation and its solutions.
It turns to be a proper extension of the Mobius group.

Special attention was paid to the nine-term recurrence relation
for the coefficients of the Taylor series solutions
of the novel symmetric form of the general Heun's equation.

We also showed that the novel Taylor series solutions are absolutely convergent inside the circle with the unit radius.
After the simple inversion of the independent variable $z\to 1/z$ one obtains also  the Laurent series solutions
which are absolutely convergent outside the circle with unit radius.

Since the four regular singular points of the symmetric form of the general Heun's equation lie
on the unit circle, one can use the new solutions for a simultaneous equal treatment of all of them.

One can hope that this new approach will simplify
the solution of the existing basic open problems in the theory of the general Heun's functions.

The novel representation will allow also development of new effective computational methods
for calculations with the Heun's functions
which at present are still a quite problematic issue.

Analogous results for the four cases of
the confluent Neun's equation will be discussed elsewhere.

\section*{Acknowledgements}

The author is grateful to Professor Sergei Yu. Slavyanov for his helpful comments,
help with the literature, and his kind encouragement, as well as
to Professors Alexander Kazakov and Artur Ishkhanyan
for their interest in the basic results of the present talk.
Special thanks to Proffesor Oleg Motygin for careful reading of manuscript and useful
remarks.

Special thanks are also to the leader of the Maple-Soft-developers, Edgardo Cheb-Terrap,
for many useful discussions during the last years
on the properties of the Maple-Heun's functions and the computational problems with them.

The author is also thankful to the leadership of BLTP, JINR, Dubna for the support and good working conditions.
This talk was also supported by the Sofia University Foundation {\it Theoretical and Computational Physics and Astrophysics}
and by the Bulgarian Agency for Nuclear Regulation, the 2014 grant.

\end {document}